# On the Variability Estimation of Lognormal Distribution Based on Sample Harmonic and Arithmetic Means

Edward Y. Ji[*] and Brian L. Ji[**]

Abstract: For the lognormal distribution, an unbiased estimator of the squared coefficient of variation is derived from the relative ratio of sample arithmetic to harmonic means. Analytical proofs and simulation results are presented.

Methods of estimating the parameters of lognormal distribution are summarized by Aitchison and Brown[1], Crow and Shimizu[2], and others. Finney's method (1941) leads to the general solutions of uniformly minimum variance unbiased estimators (UMVUE)[3], which are expressed by the infinite series or generalized hypergeometric functions[2]. On the other hand, for the harmonic mean estimation method in numerous engineering applications, some investigators preferred other estimators to UMVUE for issues such as the negative values in the series and computational simplicity[4]. Here we show an unbiased estimator of the squared coefficient of variation based on the relative ratio of sample arithmetic to harmonic means. This work is inspired by the new semiconductor integrated circuit architecture where the harmonic and arithmetic means are collectively measured[5]. While our results were originally obtained first by the empirical fittings to Monte Carlo simulations, here we start with the analytical proofs and followed by the simulation results.

The two-parameter log-normal (LN) distribution of a positive random variable $X$ is denoted as $X \sim LN(\mu_Y, \sigma_Y^2)$, with the probability density function,

$$f_{LN}(x; \mu_Y, \sigma_Y^2) = \frac{1}{x\sigma_Y\sqrt{2\pi}} exp\left(-\frac{(ln(x)-\mu_Y)^2}{2\sigma_Y^2}\right), x > 0$$

[*] Mount Kisco, NY 10549. [**] College of Nanoscale Science and Engineering, SUNY Polytechnic Institute, Albany, NY 12203. Email: bji@sunycnse.com

where $\mu_Y$ and $\sigma_Y^2$ are the mean and variance in the logarithmic space, $Y = ln(X)$. In this paper, the population parameters of $X$ in the real space are denoted as: the arithmetic mean $\alpha \equiv E[X]$, the variance $\beta^2 \equiv E[(X-\alpha)^2]$, the coefficient of variation $C_v \equiv \beta/\alpha$, and the harmonic mean $h \equiv 1/E\left[\frac{1}{X}\right]$. Furthermore, define a ratio $\omega \equiv \alpha/h$, and a relative ratio $k \equiv \omega - 1 \equiv (\alpha/h) - 1$. For the lognormal distribution, it is known that[1] $\alpha = \exp\left(\mu_Y + \frac{1}{2}\sigma_Y^2\right)$, $h = \exp\left(\mu_Y - \frac{1}{2}\sigma_Y^2\right)$, $\beta^2 = exp(2\mu_Y + \sigma_Y^2)(exp(\sigma_Y^2) - 1)$, $C_v^2 = \exp(\sigma_Y^2) - 1$, so that $\omega = exp(\sigma_Y^2)$ and $k = exp(\sigma_Y^2) - 1 = C_v^2 \geq 0$.

Let $X_1, X_2, \ldots, X_n$, $n \geq 2$, be independent and identically distributed (i.i.d.) random variables having $LN(\mu_Y, \sigma_Y^2)$. The sample arithmetic mean is $A_n \equiv \bar{X} \equiv (1/n)\sum_{i=1}^{n} X_i$. The sample harmonic mean is $H_n \equiv \overline{X_H} \equiv n/\sum_{i=1}^{n}(1/X_i)$. The sample relative ratio is defined as $K_n \equiv \frac{A_n}{H_n} - 1$.

**Proposition 1.** For $LN(\mu_Y, \sigma_Y^2)$, the expected value of the sample relative ratio $K_n$ is,

$$E(K_n) = \frac{n-1}{n}k = \frac{n-1}{n}C_v^2 \qquad (1)$$

Proof. From the multiplicative properties[1,2], if $X_a \sim LN(\mu_a, \sigma_a^2)$ and $X_b \sim LN(\mu_b, \sigma_b^2)$ are independent, then their product is a lognormal variable, $X_a X_b \sim LN(\mu_a + \mu_b, \sigma_a^2 + \sigma_b^2)$. Since $1/X_i \sim LN(-\mu_Y, \sigma_Y^2)$, $X_1/X_2 \sim LN(0, 2\sigma_Y^2)$, so that $E(X_1/X_2) = exp(\sigma_Y^2) = k + 1$. Thus,

$$E(K_n) = E\left(\frac{1}{n^2}\sum_{i=1}^{n} X_i \sum_{j=1}^{n}(1/X_j) - 1\right) = \frac{1}{n^2}E\left(\sum_{i=j=1}^{n} 1 + \sum_{i\neq j}(X_i/X_j)\right) - 1 = \frac{1}{n^2}\left(n + n(n-1)E\left(\frac{X_1}{X_2}\right)\right) - 1$$

$$= \frac{1}{n^2}(n + n(n-1)(k+1)) - 1 = \frac{n-1}{n}k = \frac{n-1}{n}C_v^2 \qquad \blacksquare$$

**Proposition 2.** Let $\hat{k}_n = \frac{n}{n-1}K_n = \frac{n}{n-1}\left(\frac{A_n}{H_n} - 1\right)$, then $\hat{k}_n$ is an unbiased estimator of $k$ and $C_v^2$, that is,

$$E(\hat{k}_n) = k = C_v^2 \qquad (2)$$

Proof. Equation (2) is straightforward from equation (1). ∎

**Proposition 3.** For $LN(\mu_Y, \sigma_Y^2)$, the variance and standard deviation (the square root of variance) of the sample relative ratio $K_n$ are, respectively,

$$Var(K_n) = \frac{2(n-1)}{n^2} k^2 \left(1 + k + \frac{k^2}{2n}\right) = \frac{2(n-1)}{n^2} C_v^4 \left(1 + C_v^2 + \frac{C_v^4}{2n}\right) \tag{3}$$

$$sd(K_n) = \frac{k}{n}\sqrt{2(n-1)\left(1 + k + \frac{k^2}{2n}\right)} = \frac{C_v^2}{n}\sqrt{2(n-1)\left(1 + C_v^2 + \frac{C_v^4}{2n}\right)} \tag{4}$$

Proof. The variance of the sum is the sum of the covariances, $Var(\sum_{i=1}^n A_i) = \sum_{i,j=1}^n Cov(A_i, A_j)$. The similar terms in the expression of the variance of $K_n$ as the sum are collected and calculated with the covariance property, $Cov(A, B) = E(AB) - E(A)E(B)$ and the lognormal properties, for example, $E\left(\frac{X_1}{X_2}\frac{X_3}{X_2}\right) = E\left(\frac{X_1 X_3}{X_2^2}\right) = E(LN(0, 6\sigma_Y^2)) = exp(3\sigma_Y^2) = \omega^3$. We have,

$$Var\left(\frac{X_1}{X_2}\right) = Cov\left(\frac{X_1}{X_2}, \frac{X_1}{X_2}\right) = E\left(\frac{X_1}{X_2}\frac{X_1}{X_2}\right) - E\left(\frac{X_1}{X_2}\right)E\left(\frac{X_1}{X_2}\right) = exp(4\sigma_Y^2) - exp(2\sigma_Y^2) = \omega^4 - \omega^2$$

$$Cov\left(\frac{X_1}{X_2}, \frac{X_2}{X_1}\right) = E\left(\frac{X_1}{X_2}\frac{X_2}{X_1}\right) - E\left(\frac{X_1}{X_2}\right)E\left(\frac{X_2}{X_1}\right) = 1 - exp(2\sigma_Y^2) = 1 - \omega^2; \quad Cov\left(\frac{X_1}{X_2}, \frac{X_3}{X_2}\right) = \omega^3 - \omega^2;$$

$$Cov\left(\frac{X_1}{X_2}, \frac{X_3}{X_1}\right) = \omega - \omega^2; \quad Cov\left(\frac{X_1}{X_2}, \frac{X_2}{X_3}\right) = \omega - \omega^2; \quad Cov\left(\frac{X_1}{X_2}, \frac{X_1}{X_3}\right) = \omega^3 - \omega^2; \quad Cov\left(\frac{X_1}{X_2}, \frac{X_3}{X_4}\right) = 0$$

For the general case of $n \geq 4$:

$$Var(K_n) = Var\left(\frac{1}{n^2}\sum_{i=1}^n X_i \sum_{j=1}^n (1/X_j)\right) = \frac{1}{n^4} Var\left(\sum_{i=j=1}^n 1 + \sum_{i \neq j}^n \frac{X_i}{X_j}\right)$$

$$= \frac{1}{n^4}\left\{n(n-1)Cov\left(\frac{X_1}{X_2}, \frac{X_1}{X_2}\right) + n(n-1)Cov\left(\frac{X_1}{X_2}, \frac{X_2}{X_1}\right) + n(n-1)(n-2)Cov\left(\frac{X_1}{X_2}, \frac{X_3}{X_2}\right) + n(n-1)(n-2)Cov\left(\frac{X_1}{X_2}, \frac{X_3}{X_1}\right)\right.$$

$$+ n(n-1)(n-2)Cov\left(\frac{X_1}{X_2}, \frac{X_2}{X_3}\right) + n(n-1)(n-2)Cov\left(\frac{X_1}{X_2}, \frac{X_1}{X_3}\right)$$

$$\left. + n(n-1)(n-2)(n-3)Cov\left(\frac{X_1}{X_2}, \frac{X_3}{X_4}\right)\right\}$$

$$= \frac{1}{n^4}\{n(n-1)(\omega^4 - \omega^2) + n(n-1)(1 - \omega^2) + n(n-1)(n-2)(\omega^3 - \omega^2) + n(n-1)(n-2)(\omega - \omega^2)$$

$$+ n(n-1)(n-2)(\omega - \omega^2) + n(n-1)(n-2)(\omega^3 - \omega^2)\}$$

$$= \frac{2(n-1)}{n^2} k^2 \left(1 + k + \frac{k^2}{2n}\right)$$

The cases of $n = 2$ and 3 are straightforward. ∎

**Example 1.** For the case of $n = 2$:

$$Var(K_2) = \frac{1}{2^4}Var\left((X_1 + X_2)\left(\frac{1}{X_1} + \frac{1}{X_2}\right)\right) = \frac{1}{2^4}Var\left(2 + \frac{X_1}{X_2} + \frac{X_2}{X_1}\right) = \frac{1}{2^4}Var\left(\frac{X_1}{X_2} + \frac{X_2}{X_1}\right)$$

$$= \frac{1}{2^4}\left\{Cov\left(\frac{X_1}{X_2},\frac{X_1}{X_2}\right) + Cov\left(\frac{X_2}{X_1},\frac{X_2}{X_1}\right) + Cov\left(\frac{X_1}{X_2},\frac{X_2}{X_1}\right) + Cov\left(\frac{X_2}{X_1},\frac{X_1}{X_2}\right)\right\} = \frac{1}{2^4}\{2(\omega^4 - \omega^2) + 2(1 - \omega^2)\}$$

$$= \frac{1}{8}(\omega^2 - 1)^2 = \frac{1}{8}k^2(k+2)^2 = \frac{1}{2}k^2\left(1 + k + \frac{k^2}{4}\right) \quad \blacksquare$$

**Proposition 4.** The variance and standard deviation of $\hat{k}_n$ are, respectively,

$$Var(\hat{k}_n) = \frac{2}{n-1}k^2\left(1 + k + \frac{k^2}{2n}\right) = \frac{2}{n-1}C_v^4\left(1 + C_v^2 + \frac{C_v^4}{2n}\right) \tag{5}$$

$$sd(\hat{k}_n) = k\sqrt{\frac{2}{n-1}\left(1 + k + \frac{k^2}{2n}\right)} = C_v^2\sqrt{\frac{2}{n-1}\left(1 + C_v^2 + \frac{C_v^4}{2n}\right)} \tag{6}$$

Proof. Equations (5) and (6) are straightforward from equations (3) and (4). ∎

Monte Carlo simulation data of $\hat{k}_n$ for various sample size $n$ and $C_v$ values are shown in Figure 1, with analytical fittings to validate equations (2) and (6). The simulations were programmed with R language[6]. The lognormal random variables are generated with *meanlog* = 0 and various *sdlog* or $C_v$ values, using the *lnorm* function. The expected value and standard deviation in Figure 1 are plotted on the same scale for guidance to the engineering application designs.

A traditional way to study the estimator efficiency is the large-sample variances[1]. Finney's UMVUE generally apply to the parameter form, $\theta_{a,b,c} = \sigma_Y^{2c}\exp(a\mu_Y + b\sigma_Y^2)$, therefore the parameter $\omega = \exp(\sigma_Y^2)$ is a special case of "$a = 0; b = 1; c = 0$". For $\hat{\omega}_{UMVUE}$, following the steps shown on page 37 of Crow and Shimizu[1], we obtain the special case solution, $Var(\hat{\omega}_{UMVUE})_{n\gg1} \to \frac{2\sigma_Y^4}{n}\exp(\sigma_Y^2)$. Since $k = \omega - 1$, $Var(\hat{k}_{UMVUE})_{n\gg1} \to \frac{2\sigma_Y^4}{n}\exp(\sigma_Y^2)$. From equation (5), $Var(\hat{k}_n)_{n\gg1} \to \frac{2}{n}\left(e^{\sigma_Y^2} - 1\right)^2\exp(\sigma_Y^2)$, thus the large-sample efficiency of $\hat{k}_n$ is,

$$eff.[\hat{k}_n] = \frac{Var(\hat{k}_{UMVUE})_{n\gg1}}{Var(\hat{k}_n)_{n\gg1}} = \frac{\sigma_Y^4}{\left(e^{\sigma_Y^2}-1\right)^2} \quad (7)$$

Figure 2 shows the large-sample efficiency of $\hat{k}_n$ according to (7).

The lognormal probability density function can be alternatively expressed in terms of $g$ and $k$,

$$f_{LN}(x; g, k) = \frac{1}{x\sqrt{2\pi ln(1+k)}} exp\left(-\frac{(ln(x/g))^2}{2ln(1+k)}\right) \quad (8)$$

where $g = exp(\mu_Y)$ is the geometric mean. The geometric, arithmetic and harmonic means are related by[7] $ah = g^2$. It is easy to see that $\hat{g}_n = \sqrt{A_n H_n}$ is a consistent estimator, that is, asymptotically unbiased when $n$ goes to infinity.

The measurement cost should be evaluated for statistical sampling. One can define the measurement cost as the number of the needed measurements. The conventional method of variability estimation requires the detailed knowledge of each replicate in the sample, thus the measurement cost number for sample size $n$ is simply $n$. Using the integrated circuits of replicating devices, some sample means may be measured collectively by various laws of physics and circuits (for example, series and parallel circuits). In these cases, the sample harmonic mean $H_n$ or arithmetic mean $A_n$ can be obtained with a measurement cost number of 1. The sample estimators $\hat{k}_n$ and $\hat{g}_n$ need a measurement cost number of 2. The $n$-to-2 measurement cost number reduction is obviously significant.

**Acknowledgement**

We thank K. Liu and M. B. Ketchen for helpful discussions about the related and broader topics.

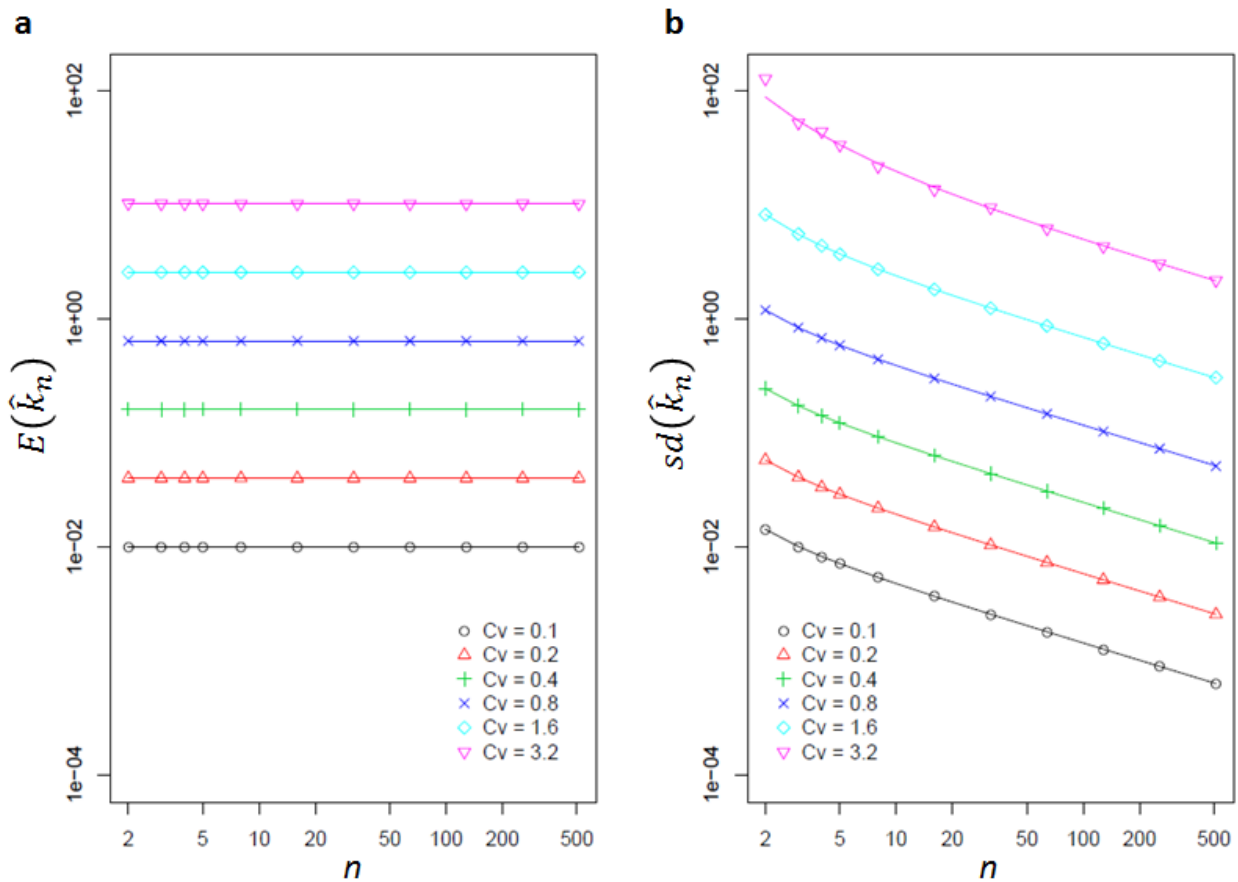

**Figure 1 | Expected value and standard deviation of $\hat{k}_n$ of various sample size *n* and various $C_v$ values.** Dots: simulation data. The expected values and the standard deviations are calculated from Monte Carlo simulation of $10^7/(n-1)$ runs. Lines: analytical predictions by equations (2) and (6).

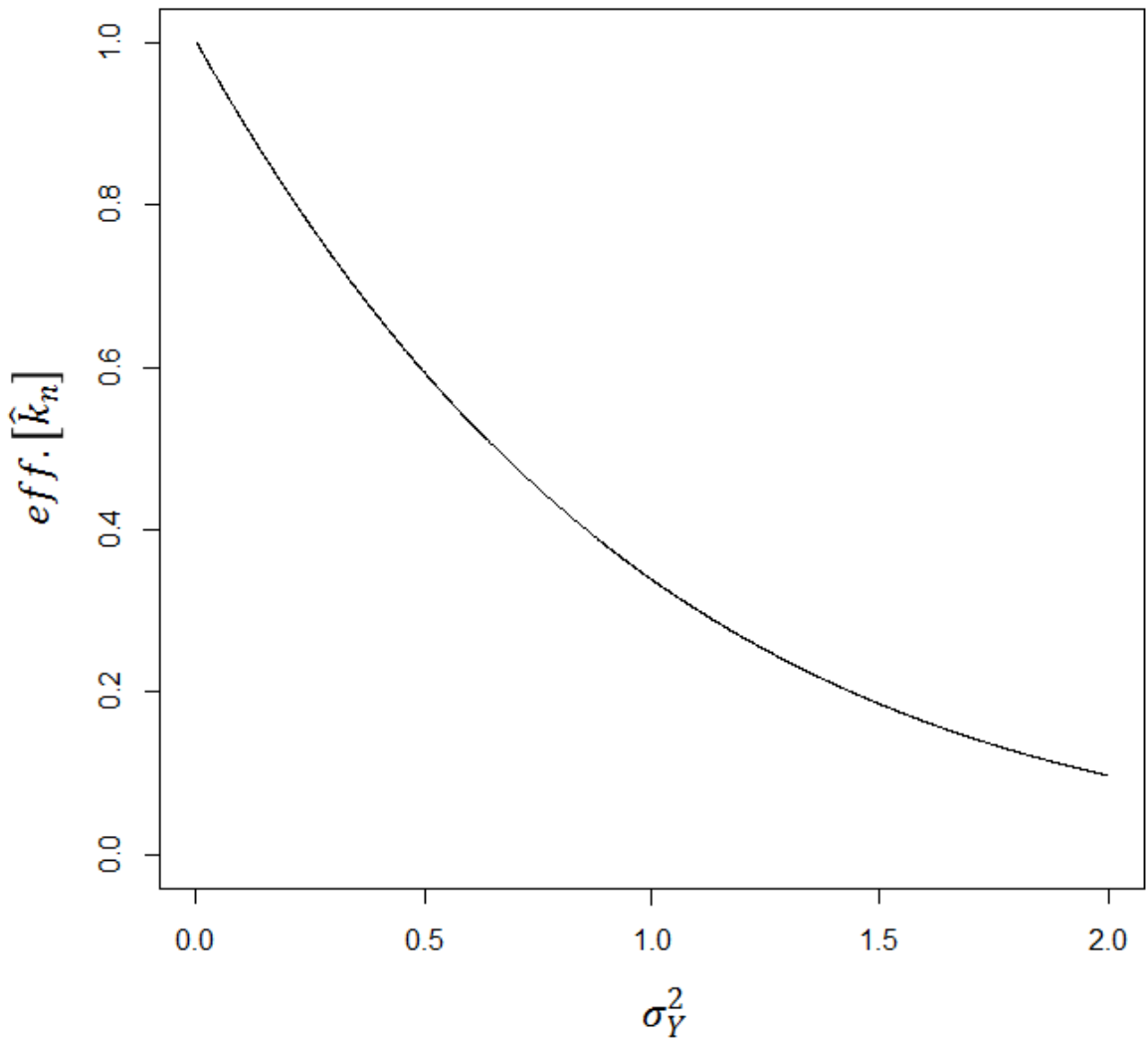

**Figure 2 | Large-sample efficiency of $\hat{k}_n$ as a function of $\sigma_Y^2$.**